*Article*

# Traditional lectures versus active learning – a false dichotomy?

**Heiko Dietrich[1] and Tanya Evans[2*]**

[1] School of Mathematics, Monash University, Australia; heiko.dietrich@monash.edu

[2] Department of Mathematics, University of Auckland, New Zealand; t.evans@auckland.ac.nz

[*] **Correspondence:** t.evans@auckland.ac.nz; Tel: +64-9-923-8783

**Abstract:** Traditional lectures are commonly understood to be a teacher-centered mode of instruction where the main aim is a provision of explanations by an educator to the students. Recent literature in higher education overwhelmingly depicts this mode of instruction as inferior compared to the desired student-centered models based on active learning techniques. First, using a four-quadrant model of educational environments, we address common confusion related to a conflation of two prevalent dichotomies by focusing on two key dimensions: (1) the extent to which students are prompted to engage actively and (2) the extent to which expert explanations are provided. Second, using a case study, we describe an evolution of tertiary mathematics education, showing how traditional instruction can still play a valuable role, provided it is suitably embedded in a student-centered course design. We support our argument by analyzing the teaching practice and learning environment in a third-year abstract algebra course through the lens of Stanislav Dehaene's theoretical framework for effective teaching and learning. The framework, comprising "four pillars of learning", is based on a state-of-the-art conception of how learning can be facilitated according to cognitive science, educational psychology, and neuroscience findings. In the case study, we illustrate how, over time, the unit design and the teaching approach have evolved into a learning environment that aligns with the four pillars of learning. We conclude that traditional lectures can and do evolve to optimize learning environments and that the erection of the dichotomy "traditional instruction versus active learning" is no longer relevant.

**Keywords:** traditional instruction; explicit instruction; explanations; active learning; undergraduate mathematics

## 1. Introduction

A recent comprehensive review of research on the teaching of proof-based undergraduate mathematics by Melhuish, Fukawa-Connelly, Dawkins, Woods and Weber [1] synthesized the findings of 104 published reports from a range of countries and research programs. Classroom pedagogy was a key differentiating characteristic chosen by the researchers to categorize the articles dichotomously. They considered two types of learning environments: (1) studies involving lecture-based pedagogy and (2) those focusing on student-centered pedagogy. The distinction was made between student-centered instruction and lecture-based instruction based on students' role in generating disciplinary ideas. In the lecture-based classroom, the instructor's primary role is to introduce mathematical knowledge; in the student-centered classroom, students (are intended to) play a substantial role in generating disciplinary ideas - they are expected to create 'new to them' mathematics to engage in creative processes akin to those experienced by professional mathematicians. It is important to clarify that this distinction differs from a much-discussed dichotomy between active and passive learning [2] which focuses on the extent to which students are prompted to engage in the learning process actively. For instance, a traditional lecture that involves quizzing, 'think-pair-share' collaborative problem-solving, or self-explanation prompts after a lecturer explained new mathematics would not be considered passive. Yet, under the categorization of Melhuish et al. [1], this type of classroom would be regarded as lecture-based and not student-centered.

To elaborate on this distinction, we utilize a four-quadrant model of educational environments conceptualized in a two-dimensional space (Figure 1). Thus, we propose a more nuanced consideration of learning environments along two axes: teacher behaviour and intended student activity. The horizontal axis represents a continuum of the teacher-facilitator role depicting the extent to which explicit explanations of disciplinary ideas are provided. Arguably, it is reasonable to consider the following clear dichotomy – either a teacher offers explanations introducing the majority of concepts and relationships between them or not. In the latter case, a teacher who believes in the value of inquiry-based learning assumes the role of a facilitator of student learning and sets up a learning environment to allow students to discover 'new to them' mathematics by themselves. The premise of this approach is that by doing so, students get to experience creative processes that are much enjoyed by professional mathematicians, thereby positively influencing their interest, motivation and, ultimately, understanding [1, 3, 4].

The vertical axis represents another dimension to account for the extent to which students are prompted to engage by the design of the learning environment. The two quadrants on the right, representing the extremes on the vertical spectrum, encompass the distinction between structured inquiry-based learning (top-right) versus unassisted discovery (bottom-right). A structured-inquiry facilitator guides learners through a sequence of specifically designed tasks to scaffold student discovery. The tasks are designed so that students may prove a major theorem or (re)invent a mathematical definition or procedure over weeks of instruction [3]. In contrast, in an unassisted inquiry-based classroom, a learner is expected to approach learning as a task of discovering something rather than learning about 'it' [5] and "to carry out his learning activities with the autonomy of self-reward or, more properly by reward that is discovery itself" ([6], p.17).

The two quadrants on the left emanate from traditional instructional approaches. A teacher plays a key role in introducing new disciplinary ideas by offering expert explanations while motivating and defining new concepts and emphasizing the relationships between newly defined concepts and the other concepts already learned. Generally, students are not expected to discover 'new to them' mathematical concepts, which is the main differentiating characteristic of the learning environments on the left-hand side in opposition to the two quadrants on the right.

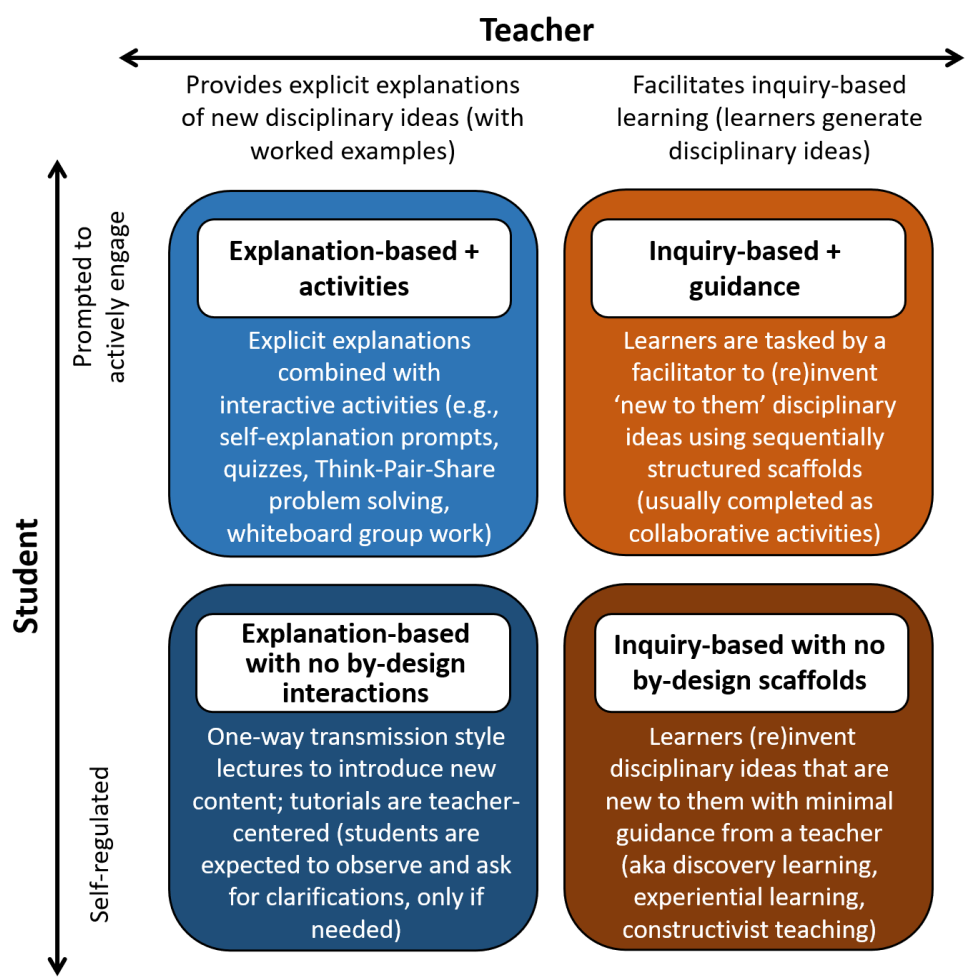

**Figure 1.** Four-quadrant model of educational environments

A serious confusion in the literature can be attributed to projecting the two-dimensional model to one dimension, in which explanation-based learning environments are viewed as belonging to the bottom-left quadrant only, counterposed to active learning environments, which combine three other quadrants simultaneously. This is evident in the large body of studies that overwhelmingly depict teacher-centered instruction as inferior compared to the desired student-centered models based on active learning techniques [2, 3]. On closer inspection, the majority of evidence reported across numerous studies points to the extra benefit of active engagement during class compared to passive listening, which is often experienced in transmission-style lectures (bottom-right quadrant). For example, the much-cited meta-analysis undertaken by Freeman et al. [2] synthesized available evidence on comparing the 'traditional lecturing' versus 'active learning' in undergraduate science, technology, engineering, and mathematics (STEM) courses. The effect sizes from 255 studies demonstrated that average examination scores were higher by about 6% in active learning sections and that students were 1.5 times more likely to fail in classes with traditional lecturing than in classes with active learning. It is important to note that the criteria for admission in the meta-analysis included equivalence in measures of performance (identical exams in both conditions), participant equivalence (randomization or quasi-randomization among treatments) and instructor equivalence (e.g., identical instructors in both treatments). These results are crucially insightful and provide strong empirical evidence against lecturing in a transmission style with anticipation of passive student engagement. However, the researchers only employed differentiation with respect to prompts for active participation (vertical dimension in Figure 1) but did not account for the difference in the provision of explanations (horizontal dimension).

Under the banner of active learning, major efforts have been made recently to promote inquiry-based mathematics learning at the tertiary level. An influential movement called Inquiry-Based Mathematical Education (IBME) has gained momentum [3], attracting substantial funding from US government agencies [4]. This resulted in the successful mobilization of regional consortia in many states, uniting over 800 mathematics education practitioners working to reform undergraduate education. IBME is characterized by the fundamental premise that learners should be allowed to learn 'new to them' mathematics without being taught, representing learning environments on the right-hand side in Figure 1. This progressive idea is based on the assumption that it is best to advance learners to the level of experts by engaging learners in mathematical practices similar to those of practicing mathematicians: creating new definitions, conjectures and proofs - that way, learners are thought to develop 'deep mathematical understanding'. Despite heavy criticism leveled at this approach from the cognitive science perspective [4, 7, 8], the IBME proponents continue to unreservedly promote their approach by referring to the evidence in support of active learning vs traditional lecturing [3]. Again, this conflation is the result of omitting to account for the two-dimensional complexity of learning environments, as illustrated by the four-quadrant model in Figure 1. We return to this point in the Discussion section.

In summary, the experimental evidence does not indicate that the top-left quadrant, an explanation-based learning environment with prompts for active engagement, is ineffective [2]. But it seems that this point is lost at present, given that, in general, researchers in undergraduate mathematics education frequently advocate that "the teaching of undergraduate mathematics courses should rely less on lecture and more on student-centered inquiry-based instruction" ([1], p. 4). This resulted in the erroneous dichotomous stance, which is the subject of our study.

Consequently, educators who still see merit in traditional instruction to provide explanations about new definitions, theorems, proofs etc., often feel pressured to justify their choice. The aim of this case study is to present one such justification. While being strong advocates for adapting teaching methods that have scientifically been shown to be effective, we argue that traditional lectures have evolved and, therefore, can have, under certain assumptions, a place in some areas, such as high-level mathematics. For this, we first provide a brief overview of recent research on effective teaching and learning, based on Stanislav Dehaene's theoretical framework, comprising "four pillars of learning" [7]. Through the lens of this theory, we then analyze a case study of how traditional instruction can successfully be

5embedded in a student-centered unit design. Tracking the changes implemented by the first author, who has been teaching abstract algebra courses using traditional instruction for more than a decade, we illustrate how over time, the unit design and teaching approach has naturally evolved into a learning environment that, coincidently, seems to align with the four pillars of learning.

## 2. How We Learn: Dehaene's Four Pillars of Learning

The French author and cognitive neuroscientist Stanislav Dehaene [7] studied extensively how learning changes the brain and how the learning process can be facilitated. In his 2020 book "How we learn", he identified four major functions that are essential for effective learning. Due to their importance, Dehaene calls them the four pillars of learning: "Teachers who manage to mobilize all four functions in their students will undoubtedly maximize the speed and efficiency with which their class can learn. Each of us should therefore learn to master them" ([7], p. 146). We briefly summarize these four pillars before discussing our case study using this theory.

**Pillar 1: Attention**

Attention refers to "all the mechanisms by which the brain selects information, amplifies it, channels it, and deepens its processing" ([7], p. 147). The human brain constantly receives a myriad of signals, and attention helps it to select the relevant information and to strengthen the activity of relevant synapses to store the new information deep in its cortical circuits; conscious attention facilitates this "long-term potentiation" ([7], p. 150). It is clear that without attention or with misguided attention, the learning process is significantly hindered [9]. Dehaene summarizes, "This is why every student should learn to pay attention – and also why teachers should pay more attention to attention!" ([7], p. 150). According to the American psychologist Michael Posner, three major attention systems are relevant: these are Alerting, Orienting, and Executive Attention [10]. The first attention system is responsible for the initial process of being alerted to a new situation; it tells the body when to increase its cognitive ability to be vigilant and focused. The second system, orienting, is responsible for selecting which input signals are receiving appropriate attention; it tells the brain what to focus on. It is therefore recommended that teachers carefully choose where they want to direct the attention because "only the items that lie at the focus of attention are represented in the brain with sufficient strength to be efficiently learned" ([7], p. 158). Items outside the focus of attention could be completely suppressed, which is well illustrated by the famous "invisible gorilla" experiment [11]. The third system, executive attention, is responsible for deciding how the information that is in the focus of attention is processed; it is said to correspond to what is called "concentration" ([7], p. 162). Dehaene compares this system with the switchboard of the brain, which "orients, directs, and governs our mental processes" ([7], p. 159). The way how the brain responds to a certain stimulus is crucial to learning because, as Dehaene declares, our brain is not able to truly multitask. Learning to concentrate is, therefore, an important part of the learning process: "Teaching requires paying attention to the limits of attention and, therefore, carefully prioritizing specific tasks." ([7], p. 162). In summary, in the context of education, learning can be facilitated if the educator manages, firstly, to captivate their students' attention by constantly stimulating their interest (alert system), then guiding the students to where exactly they should focus their attention on (orientation system), while constantly considering what the students do not know to assist their executive systems in making the most appropriate choices (executive attention).

**2.1. Pillar 2: Active engagement**

Dehaene describes that external stimuli alone are not sufficient for learning but that active engagement is required: "efficient learning means refusing passivity, engaging, exploring, and actively generating hypotheses and testing them" ([7], p. 179). Importantly, active engagement takes place solely in our brain while it is focused and actively generates new mental models. This means that students who watch a lecture might seem passive from the outside, but their brains could be alert and active, working hard to internalize new concepts by constantly rephrasing them and by updating mental models of the world. Even though Dehaene references studies that demonstrate the inferiority of traditional lecturing in STEM, where students remain passive while the teacher "preaches" for fifty minutes ([2] as cited in [7]), he emphasizes that this is not to say that a traditional instruction intrinsically prevents active engagement. He suggests that "[t]here is no single miraculous method, but rather a whole range of approaches that force students to think for themselves" ([7], p. 181). Dehaene lists curiosity, the desire to learn, as one of the key ingredients to active engagement: as the "brain's governor", it guides us to what we think we can learn by identifying gaps between what we already know and what we would like to know ([7], p. 190). To help students maintain a high level of curiosity, Dehaene suggests educators choose the appropriate level of cognitive stimulation, rewarding curiosity, and leaving room for students to self-explore their own questions ([7], p. 194-196).

**2.2. Pillar 3: Error feedback**

Studies suggest that learning is facilitated if the brain detects a gap between its expectations and observations: "surprise is one of the fundamental drivers of learning" ([7], p. 202). In this situation, error signals flag a discrepancy between the expectation and observation, even though this discrepancy can also be in the form of a *correct* expectation but with low confidence. For instance, if a student has little confidence in their solution, they might still be *surprised* when their answer is correct, and the resulting "error signal" still leads to a brain stimulus that enhances learning. "Receiving feedback improves memory even when the right answer was chosen" ([7], p.220) because "as long as knowledge is not perfectly consolidated, the brain continues to learn, even weakly." ([7], p. 220). Dehaene concludes it is important for educators to provide timely and precise feedback to their students: this will lead to surprise moments (activating the learner's brain) and will considerably enrich the information available for self-correction [7]. Dehaene stresses that "error feedback should not be confused with punishment" ([7], p. 210), which, in an educational setting, is often associated with awarding low marks; he mentions that consistently receiving low marks can lead to school-induced stress, and "mathematics anxiety" is a particular example of that ([7], p. 212). It is therefore suggested that testing should not primarily be seen as a tool for grading students but rather as a tool for implementing a retrieval practice, which is deemed "one of the most effective educational strategies" [12-14]. Testing one's knowledge will automatically lead to many of the aforementioned surprise moments, which, in turn, facilitates learning. Moreover, these retrieval episodes should be spaced out according to the benefits of the distributed practice (spacing) effect, which is a well-documented effect in the experimental psychology literature [15, 16]. Dehaene recommends that testing and learning to be spaced out over a more extended period to ensure retention in the long term: "…each review reinforces learning. It refreshes the strength of mental representations and helps fight the exponential forgetfulness that characterizes our memories" ([7], p. 219).

### 2.3. Pillar 4: Consolidation

The last pillar of learning, consolidation, describes the "shift from slow, conscious, and effortful processing to fast, unconscious, and automatic expertise" ([7], p. 222). Examples include the ability to read texts fluently and effortlessly. Consolidation goes hand in hand with repetition and spacing out of learning, which is also connected to the nightly breaks during which we sleep: Sleep is when our brain "replays the important events it recorded during the previous day and gradually transfers them into a more efficient compartment of our memory" ([7], p. 224). During sleep, "our floating ideas from the day are reactivated hundreds of times at an accelerated rate, thus multiplying our chances that our cortex eventually discovers a rule that makes sense" ([7], p. 231). This supports anecdotal stories in which learners, who were completely immersed in a topic during the day, *made discoveries* during their sleep.

## 3. A Case Study: Traditional Lectures and Active Learning

We now present a case study of the evolution of the learning environment created by the first author (Heiko), who has been teaching abstract algebra courses (also known as units) using traditional lectures for a decade. By a "traditional lecture" we mean a form of instruction where a lecturer's main focus is on providing expert explanations of new material to students with only limited time spent on student-lecturer interactions. We argue that even with this traditional approach, a suitable unit design allows Heiko to run an effective and student-centered class. We illustrate this by reporting on Heiko's experience with the third-year algebra unit MTH3150 that he has been fully responsible for at Monash University since 2013; the continuity of teaching the same class to the same level of students is useful for analyzing under the case study methodology when describing the impact and value of the approaches outlined below.

### 3.1. Background information about MTH3150

MTH3150 (formerly "Algebra and Number Theory 2", now "Algebra 2: Rings and Fields") is a Level 3 pure mathematics unit at Monash University. Most of the students are in their third year of completing a major in mathematics. The content of MTH3150 comprises basic concepts in ring theory, such as integral domains, principal ideal domains, unique factorization domains, Euclidean domains, field extensions, and applications in coding theory. The focus is on understanding abstract mathematical concepts and proofs, making MTH3150 one of the more abstract units in Monash's mathematics undergraduate programme. This poses various challenges. First, students are often overwhelmed by the level of mathematical abstraction and the amount of time spent on mathematical proofs. Second, the pace of this unit is fast, with almost every lecture introducing new definitions and concepts.

The workload of MTH3150 is distributed into three weekly lectures (3 hours), one weekly applied class (90-120 min), and time for personal study (7 hours for revision and written assignments). Students have access to a unit website (where they can ask questions in a forum) and weekly consultation hours. During the lectures, Heiko presents the lecture content in a traditional style of instruction. The applied classes run in a format known as "whiteboard tutorials": students gather in small groups (usually 3-4 students) in front of a whiteboard and work collaboratively on a pre-designed problem sheet. There are also weekly written assignments whose primary purpose is not merely an assessment of learning but facilitating and guiding timely and regular revisions, as well as consolidation of learning. It is important to note that the three lectures only constitute one quarter (3 hours) of the nominated workload (12 hours per week) for MTH3150; the other parts of the unit contain multiple active-learning components.

Next, we present an analysis of the evolved learning environment in this unit through the lens of Dehaene's teaching and learning framework outlined above.

**3.2. Analysis of learning environment: Pillar 1 (attention)**

During the lecture component of the unit, several strategies are used to raise students' attention. This starts before each lecture commences because the first slide (shown while Heiko waits for the lecture to start) usually contains a topical mathematical riddle or cartoon. Students become visibly engaged when this slide is up, which often leads to an engaged pre-lecture small talk that effortlessly transitions into the first part of the lecture (usually a recap section). Recalling what has happened in the previous lectures (including prompting students) seems an effective way to raise attention at the beginning of a lecture. Each lecture is delivered via handwriting on a tablet computer and accompanied by many verbal explanations, alternating auditory/verbal and visual channels to facilitate the effectiveness of attention in sensory memory. This is explained by Clark and Paivio's dual-coding theory [17] and Baddeley and Hitch's model of working memory [18]. Working memory has two systems: verbal/auditory processing takes place in the so-called phonological loop, whereas visual processing takes place within the visuospatial sketchpad (e.g., [19]). The two systems are complementary, meaning that if processing can be split between the two systems, then the total working memory capacity could be increased. This is the mechanism utilized by Heiko in order to reduce the processing load of the working memory of the students when introducing new material.

Handwriting slows down the pace of the lecture and makes Heiko more attentive: when writing down a mathematical argument by hand, he usually repeats the thought process that led to it and shares this experience with the students. Heiko also encourages his students to take notes because the latter is still considered an important aspect of formal classroom learning associated with higher achievement [20]. Handwriting also gives Heiko the flexibility to make quick amendments or respond to questions during a lecture, allowing for flexibility in interactions during the lecture.

Heiko uses highlighters and color coding for his writing (for example, definitions and formulas are blue, exercises are red, etc.), which effectively draws learners' attention to particular parts of the writing. Heiko often focuses students' attention on key features of concepts or parts of proofs by asking quick (and sometimes rhetorical) questions: "What is wrong here?", "Why is that?", "Why is this important?"; sometimes, individual students are asked directly, "Xuan, do you have an idea of what the next argument could be?". Body language, moreover, is a powerful tool to raise attention: Stepping away from the lectern, underpinning words with hand gestures, and maintaining a lively classroom persona are important tools for Heiko to keep students attentive. This includes varying his intonation, repeating and reinforcing statements, including theatrical pauses for raising suspension and maintaining regular eye contact with the class. This is particularly effective when drawing students' attention to a surprising result or argument. At the end of some proofs, Heiko often steps back from the lectern and provides a conceptual overview of the main steps used in the argument; this is another example of how Heiko directly utilizes the attention systems "orientation" and "executive attention".

In summary, despite the lack of established active-learning exercises during the lectures (such as Think-Pair-Share moments, group discussions, or peer review), Heiko utilizes a variety of techniques to alert students and to direct their thought processes, which, according to the cognitive science research is beneficial in enabling effective learning.

### 3.3. Analysis of learning environment: Pillar 2 (active engagement)

Having raised the attention of the students in class, Heiko stimulates their curiosity by asking questions and making students think about the problems and concepts introduced during the lesson. For example, he often rephrases his writing and then asks students to repeat the argument in their own words. If a follow-up explanation is required, then Heiko can quickly provide this (see the comment on handwriting above). Even though not much time is spent on genuine discussions, these small interludes often give the lecture an interactive feel: students know they can ask questions, and the lecturer will respond. In addition to that, Heiko often includes additional exercises in his notes. For example, there might be a statement for finite fields, and the exercise asks: "What goes wrong for infinite fields?". These questions are meant to be additional stimuli for the students to think about during their revisions (not necessarily during class). Of a similar nature is the additional materials that Heiko provides via the unit website: supplementary resources and additional questions are uploaded, encouraging student participation in the online discussion forums.

However, the parts of the unit that are structured around active (behavioral and cognitive) engagement are the applied classes (also known as tutorials) and weekly assignments. At the beginning of each applied class, students are asked to discuss with each other the terms and definitions required for each question on the problem sheet; this serves as a quick revision to ensure everyone understands the context of the question. Subsequently, the students work on the questions in group-research-style by developing a solution to a given problem together on the whiteboard. The tutor is often merely an observer-facilitator; they walk around from group to group, give hints, make corrections, and may ask follow-up questions to groups that have progressed particularly well. These problem-solving sessions are designed to provide opportunities for learners' active (cognitive) engagement, which according to cognitive science research, is beneficial for generative learning [21]. They happen *after* the relevant material has been discussed in class, which seems well justified according to a study by Ashman, Kalyuga and Sweller [22] that concluded that in learning environments where "element interactivity" (defined as the complexity of a concept in terms of its connections to and dependence on other concepts) is high, it is beneficial to *first* have explicit instruction and *then* problem-solving sessions.

MTH3150 has weekly assignments where students are asked to submit their work on specific problems. Students are encouraged to discuss the problems with their peers, but collusion and plagiarism are not allowed. These weekly problem sheets motivate students to keep up with their revisions, thereby continuously enabling them to stay actively engaged with the unit.

Since the pace of the lectures is high, there is often insufficient time during class for students to completely understand all new concepts. For some students, this is a new experience that can be detrimental to their motivation. To increase motivation and thereby potential active engagement, Heiko communicates the expectation that understanding each lecture requires a thorough revision (that will be guided via applied classes and assignments). Such explicit alerting is designed to provide reassurance and act as an incentive to revise the content.

**3.4. Analysis of learning environment: Pillar 3 (error feedback)**

During the lectures, Heiko often creates "surprise" moments when he asks a question and then reveals solutions (sometimes purposely incorrect solutions to make students think). This is in line with the third pillar of learning put forward by Dehaene, who emphasized that learning is facilitated if the brain detects a gap between its expectations and observations; he posited that surprise is one of the key drivers of learning [7]. As mentioned in Sections 3.2 and 3.3, asking quick questions and making students think about their answers is common during Heiko's lectures and applied classes. During the applied classes, students also receive continuous feedback on their thinking process and writing from their working group and the tutor. Each written assignment, moreover, is marked by a tutor and returned to the student with personalized written feedback. Since the assignments in MTH3150 are weekly, each worth only about 4% of the final unit mark, this shifts the role of the assignments from assessment tools to an effective tool for enhancing the learning experience. Regular assignments motivate students to revise the lecture content in a timely manner, and it allows students to get frequent feedback on their learning.

Self-testing is also encouraged: In addition to recommending one of the most effective learning strategies - retrieval practice [13] (e.g., by using flash cards) - Heiko releases lists with "Open Review

Questions". These questions are open-ended questions that eventually cover all topics of the unit and are designed as prompts for self-explanations. A large body of experimental studies demonstrated the generative learning benefits derived from engaging in the self-explanation practice, whereby students are prompted to explain to themselves why a particular logical step is true [23]. Furthermore, Heiko utilizes the benefits of another research-based effective strategy - elaborative interrogation [14] - by asking students to elaborate on what they already know about a specific concept. For example, in the Open Review Questions, he asks students to think about integral domains: they are prompted to state the formal definition, give examples *and* non-examples, describe basic properties, and think about situations where they have used these concepts in class. Some students even use these questions during applied classes to test each other's knowledge.

**3.5. Analysis of learning environment: Pillar 4 (consolidation)**

The consolidation process is supported in MTH3150 by the following means. First, learning activities are spaced out over multiple weeks. If a topic is introduced during the lectures in Week *n*, then this topic is discussed (directly or indirectly) in applied classes and written assignments during Weeks *n*, *n+1*, and *n+2*. Second, many concepts in MTH3150 are built upon each other (for example, the discussion of polynomial rings employs principal ideal domains, fields, maximal ideals, irreducible elements, etc.), which enables distributed (spaced) practice. As mentioned previously, the distributed practice effect is a well-evidenced effect pertaining to improved learning through effective consolidation when learning episodes are separated by non-study intervals – in other words, it is the opposite to *cramming* (for review, see [16]). The entire structure of the MTH3150 unit designed by Heiko, which includes spaced repetition during lectures, followed by prompts to retrieve and revise during the applied classes (whiteboard tutorials), weekly assignments and Open Review Questions, assists students with this process.

**3.6. Student perspectives and outcomes**

The structure of the unit designed and implemented by Heiko is appreciated by the students. To give an indication of the student perspectives, we summarize qualitative and quantitative student evaluation results from 2013-20. The upper part of Table 1 summarizes some of the quantitative evaluations (each score is the median of the responses, each out of 5.00), providing data for the following items:

- **Activities** refers to "*The tutorials helped me achieve the learning outcomes*" (2013-15) and "*The Activities helped me achieve the learning outcomes for the unit*" (2016-20).
- **Assessment** refers to "*The assessment tasks helped me achieve the learning outcomes*" (2013-15) and "*The assessment in this unit allowed me to demonstrate the learning outcomes*" (2016-20).
- **Unit** refers to "*Overall I was satisfied with this unit*" (2013-20).
- **Teaching** refers to "*Overall I would rate Heiko's teaching*" (2013-15) and "*In relation to MTH3150,*

*I found Heiko clear in his explanations*" (2016-18) and "*Overall I was satisfied with Heiko's teaching*" (2019-20). The row **"% top"** relates to **Teaching** and lists the percentages of responses that were "*good/outstanding*" (2013-15) and "*(strongly) agree*" (2016-20).

The average SET response rate for MTH3150 between 2013-20 is 57.8%, which is well above the Faculty of Science average (usually between 42-50%). However, we recognize a limitation of the study in reporting only the data from students who have chosen to provide an evaluation. Whenever available, the average Faculty scores are also listed; due to different reporting forms, this information was not available for all years. Student appreciation is reflected in the high scores for all evaluation items, with the highest scores reported for **Activities** and **Teaching**, which provides some evidence in support of our claim that traditional unit structures can evolve to accommodate students' needs. The gradual improvement of the unit design and teaching practice is evident in the variation of scores over the eight-year period, with the absolute maximum in teaching scores reached in the last two years of the period.

**Table 1.**
*Student Evaluations and Outcomes for MTH3150*

| Year | 2013 | 2014 | 2015 | 2016 | 2017 | 2018 | 2019 | 2020 |
|---|---|---|---|---|---|---|---|---|
| **Response** | | | | | | | | |
| # students | 35 | 30 | 21 | 46 | 20 | 28 | 32 | 35 |
| # responses | 17 | 17 | 12 | 27 | 14 | 17 | 20 | 17 |
| Response rate | 48.6% | 56.7% | 57.1% | 58.7% | 70% | 60.7% | 62.5% | 48.6% |
| **Unit Evaluation** | | | | | | | | |
| Activities | 4.43 | 4.65 | 4.58 | 4.54 | 4.64 | 4.65 | 4.59 | 4.85 |
| Assessment | 4.56 | 4.79 | 4.50 | 4.60 | 4.57 | 4.79 | 4.79 | 4.83 |
| Unit | 4.44 | 4.56 | 4.64 | 4.32 | 4.50 | 4.56 | 4.83 | 4.79 |
| *Faculty average* | *4.03* | *4.03* | *n/a* | *4.06* | *4.08* | *4.17* | *n/a* | *4.15* |
| **Teaching Evaluation** | | | | | | | | |
| Teaching | 4.70 | 4.65 | 4.83 | 4.58 | 4.42 | 4.80 | 5.00 | 5.00 |
| *Faculty average* | *n/a* | *n/a* | *n/a* | *4.34* | *4.35* | *4.45* | *n/a* | *4.51* |
| % top | 81.3% | 94.1% | 100% | 95.8% | 100% | 100% | 100% | 100% |
| **Outcomes** | | | | | | | | |
| # students with exam | 31 | 28 | 18 | 39 | 19 | 26 | 29 | 31 |
| HD | 22.9% | 42.9% | 27.8% | 17.9% | 26.3% | 26.9% | 37.9% | 38.7% |
| D | 19.6% | 10.7% | 16.7% | 30.8% | 10.5% | 30.8% | 13.8% | 19.3% |
| N | 16.1% | 21.4% | 22.2% | 12.8% | 21.1% | 11.5% | 6.9% | 12.9% |
| Mean (final unit mark) | 62.6% | 67.7% | 63.8% | 63.9% | 61.4% | 69.0% | 70.0% | 71.5% |

*Note*. HD = High-Distinction (grade band: 80-100%), D = Distinction (70-79%), N = Non-pass ( < 50%)

The lower part of Table 1 comments on the unit outcomes: it lists the percentages of the students that have written the final exam and reached a high-distinction (HD, 80-100%), a distinction (D, 70-79%), and a non-pass (N, <50%) as their final unit grade; the last row lists the mean of the final unit mark. (Students that have not sat the final exam usually have withdrawn from the unit or they had special consideration and were eligible for a deferred examination paper at a later time.) The variation in the fail rate is explained by the fact that the class size is relatively small, and each year there are a few students that do not submit all assessment tasks (and consequently lose valuable marks). Those students who engage with the unit, submit all assessment tasks, and attend all applied classes, usually pass the unit.

Further evidence demonstrating how a traditional mathematics unit can evolve to become more student-centered is provided in students' open-ended responses as part of the MTH3150 formal evaluation. The following anonymous student quotes are from responses to the questions "*What were the best aspects of the unit?"* and "*What aspect(s) of this unit did you find most effective*?", respectively.

**Lectures:** *I liked that [Heiko] did a recap of the previous lecture every time too, it was a good refresher (...) Heiko's cartoons, mini math problems and quizzes at the beginning of the lectures were a nice and interesting way to begin each lecture (...) [Lectures] were great, they were stimulating and they were well presented (...) [Heiko] interacts with students in the lectures which is excellent (...) The lectures and tutorials are engaging (...) The fact [Heiko] hand writes the notes [...] made learning so much easier.*

**Assessment:** *The structure of weekly assignments reinforcing the content learned in the previous week I found very effective as it ensured we all stayed on top of the content all the times (...) The weekly assignments, as they provided a way for me to see how well I was progressing in the unit (...) The weekly assignments are perfect for consolidating knowledge throughout the semester. Please do not get rid of them! (...) I really enjoy having assignments due each week. As much as it is stressful at the time, I think it is invaluable in terms of my learning throughout the semester. The feedback given on each assignment is also extremely helpful.*

**Applied classes:** *I found the tutorials an excellent resource/opportunity to challenge my understanding of the content (...) The tutorials are the most engaging and educational I have experienced in my maths degree. It ensured we worked and learnt together and were a pleasure to be a part of. (...) I thought that the tutorials were great as they were interactive, and were much more effective than the traditional maths tutorials (...) Working on the whiteboard in groups during the tutorials is by far the best way to learn the material (...) Tute style is the best of any of the maths units I've done (...) I also genuinely enjoyed the tutorials which is a first I can say for a Maths unit (...) The tutorials were brilliant - I found the group setting and whiteboard work was really effective to consolidate what I'd learned and also find out where I wasn't understanding a concept.*

These comments provide some support to the hypothesis that Heiko's changes prompting students to actively engage more often have improved learning environment. However, it is important to note that care should be taken in interpreting the results of Student Evaluation of Teaching (SET) since recent studies (e.g., Uttl et al. [24]) indicate that there is often no or only minimal correlation between SET ratings and learning; indeed, outcomes are sometimes influenced by curriculum-unrelated interventions (such as distribution of cookies, see Hessler et al. [25]). Nevertheless, the consistency of the data presented above indicates, to a certain extent, that the learning environment and teaching model of the unit seem well appreciated by the students. Moreover, the longitudinal data on learning outcomes demonstrates that as Heiko worked on incremental changes in his traditional classroom, the learning outcomes have not diminished. In fact, the proportion of high-achieving students has reached its absolute maximum, with 58% of students achieving distinction or high-distinction in the last year of the studied period. As for the group of students who failed the course, the numbers overall indicate a downward trend, with a couple of outliers (2017 and 2020), which can be explained by specific cases outside of students' control and an e-exam in the first year of the COVID-19 pandemic.

After analyzing student data over the years, Heiko concluded that success in this course is predicated on continuous engagement, and if one falls behind, then it is very difficult to catch up. To address this issue specifically, this year Heiko introduced regular catch-up sessions (meetings where students can ask anything that requires clarification) and weekly "self-reflections" sessions, where students are asked to reflect on each week's lectures and to comment on what aspects were most challenging (the 'muddiest point') and what was most thought-provoking. The latter is an additional tool to monitor student progress continuously and to provide personalized responses and feedback. It is hoped that these additional innovations will help all students to reach their full potential. However, the data on the impact of these interventions is not ready to be analyzed, and, thus not reported in this study.

## 4. Discussion

Why did Heiko choose traditional explanation-based instruction for his unit? The MTH3150 content is abstract, and students usually require more time than is available during explanation-based classes (lectures) to understand new concepts fully. By this, we mean to select the new information (by paying attention) and then organize and integrate it with the existing knowledge activated from a learner's long-term memory [26]). This makes it challenging to have deep ad-hoc discussions and problem-solving activities during the lectures. The knowledge and experience of the instructor help to make sure the content of each lecture is carefully chosen and attention is drawn to the key principles and characterizing features of the definitions, concepts and arguments. This level of guidance is usually not possible when students are asked to learn the material themselves (e.g., in an Inquiry-Based classroom [4]). Indeed, it seems that the additional verbal explanations accompanying the visually presented mathematics are most valued by students, who often ask for more of these explanations to clarify intricate steps in a complicated argument.

Time constraints are another reason for the choice of traditional instruction. Without reducing the content of the unit, there is not enough time to include most of the established active learning methods in a 50-minutes lecture; this is even more difficult since the aim is that each class contains a recap, motivation, and outlook. Fundamentally, Heiko designs lectures to serve as the initial stimulus that would spark the subsequent active cognitive engagement and consolidation processes; he uses the three hours of lectures per week to provide high-quality explanations that form the basis for the active learning activities that follow. As outlined in the Introduction, there is plenty of evidence about the benefits of active learning, which is the foundation in the design of applied classes and weekly assignments in the unit. Overall, MTH3150 students spend three-quarters of the weekly allocated workload on active learning tasks.

*It is crucial* to emphasize that the goal of this study is *not* to promote a traditional course structure in which students are expected to passively listen to the lecturer's and tutor's explanations with no expectations to actively engage with mathematics except by themselves at home (bottom-left quadrant in Figure 1). Instead, we argue that the traditional learning environment in MTH3150 has evolved, over time, to include traditional lectures that actively engage students, not by employing time-consuming active learning activities but by a student-focused educator providing expert explanations with frequent prompts for students to engage in self-explanations. After all, according to research, "active engagement takes place in the brain" ([7], p. 178), and our case study illustrates that this can be stimulated during traditional instruction: the "monologue lecture" has evolved into a lecture that actively addresses and involves the audience.

Simultaneously, the format of the applied classes has changed from traditional instruction (with a tutor demonstrating how to solve problems on a board) to active-learning whiteboard tutorials (with students working collaboratively in small groups). Overall, these developments led to a new unit design that now comprises 75% of active learning activities and 25% of traditional explicit instruction, corresponding to the top-left quadrant in Figure 1. Arguably, further reduction of the explicit explanations component would be counterproductive. This is supported by the recent research examining mathematics educators' varied quality (and the lack) of explanations and their effects.

In a 2016 study, Lachner and Nückles [27] concluded that the explanations provided by mathematicians (with lower pedagogical content knowledge but high content knowledge) and mathematics teachers (with high pedagogical content knowledge but lower content knowledge) about an extremum problem intended for high school students, were principally different. The teachers mainly focused on the solution steps as an algorithm for finding extreme values of a function (product-orientation). In contrast, while demonstrating the solution steps, the mathematicians clarified why a specific step in the solution was necessary (process-orientation). In their follow-up experimental study, the researchers investigated the effectiveness of these explanations on eighty high-school students who were randomly split into three groups, receiving: (1) product-oriented explanations, (2) process-oriented explanations, or (3) no explanations (inquiry-based learning group). They found that students who were not provided with an explanation showed the lowest learning gains. Furthermore, students who learned with a process-oriented explanation performed significantly better than students who learned with a product-oriented explanation on an application test. This finding was replicated and extended in a 2019 study [28] in which the impact of different explanations was assessed in a randomized sample of 129 students receiving either principle-oriented or procedure-oriented explanations on four mathematical topics. Students provided with principle-oriented explanations substantially outperformed those given procedure-oriented explanations on the application test (with similar problems) and, importantly, on the transfer test (with dissimilar problems). This research underscores the critical role quality explanations play in learning processes and support our assumption that an optimally designed mathematics unit should contain a substantial explicit instruction component.

Despite the apparent scarcity of research identifying causal relationships between pedagogy and learning outcomes at the tertiary level [1], our perspective is in line with some tertiary mathematics literature, albeit limited in volume. For example, one of the esteemed researchers with extensive experience in teaching advanced mathematics, Lara Alcock, advocated for 'tilting the classroom' to make mathematics lectures more engaging without the need for a wholesale classroom restructure [29]. Based on two decades of rigorous empirical research, Alcock offers a suite of short-and-snappy activities that can be embedded into a traditional lecture to capture attention and foster a positive atmosphere. Recently, two other prominent researchers, Keith Weber and Timothy Fukawa-Connelly, have presented an in-depth analysis of views expressed by pure mathematicians who were learners themselves while attending a 2-week workshop on inner model theory, a branch of set theory, that was delivered as a series of lectures [30]. The attendees were convinced that the lectures were invaluable to them because they provided a roadmap outlining the central ideas of the theory, helped direct their attention to what was important and explained the utility of the technical tools used to develop the theory. Moreover, the lectures motivated them to want to understand the material, and a high-level overview served as a structure guiding further study to unpack logical details. The authors concluded that further studies are warranted to investigate the function of lectures in enabling effective learning. Moreover, some theoretical progress has recently been made by Ofer Marmur [31], who proposed the construct of *Key Memorable Events* as a theoretical and methodological lens to examine the interplay between affective states experienced by learners during traditional lectures and their relation to cognition.

As mentioned before, this is in sharp contrast with the majority of the literature comprising the field of undergraduate mathematics education, which "generally advocates that the teaching of undergraduate mathematics courses should rely less on lecture and more on student-centered inquiry-based instruction" ([1], p. 4). Specifically focusing on the nuanced distinction between inquiry-based learning and the umbrella term 'active learning' (horizontal and vertical axes in Figure 1), we have recently countered this view in a scoping review of the literature [4]. Through consolidating the evidence available from 60 years of experimental research and theoretical accounts explaining how efficiency in learning is best achieved, we concluded that the call for reform from the Inquiry-Based Mathematics Education (IBME) advocates is not justified [4]. Specifically, we assert that evidence does not support the general claim that students learn better (and acquire superior conceptual understanding) if they are not lectured. Neither is the general claim about the merits of IBME for addressing equity issues in mathematics classrooms [4].

The case study reported in this paper is not an isolated example. We observe in our professional network that whiteboard tutorials seem to be the new standard for applied classes and that many of our colleagues have adapted a similar approach to lecturing mathematics: they use traditional explicit instruction to provide explanations, but the lectures are not implemented as transmission-style monologues anymore. To illustrate, a 2019 formal teaching observation undertaken by an education leader at Monash University of an MTH3150 lecture concluded, "I do not think the lecture was overly traditional. It was actually quite active and blended" – despite the fact that most of the lecture time was spent on instruction.

We conclude with an observation that traditional lectures can evolve (and have evolved in many places) and that traditional instruction can be effectively embedded in a student-centered unit design that is research-informed and evidence-based, according to the latest findings from cognitive science, educational psychology, and neuroscience. It is, therefore, misguided to erect a dichotomy between traditional lectures and active learning, which unavoidably impedes a nuanced consideration of the full spectrum of possibilities and skews research agendas.

This research has received approval from Monash University Human Research Ethics Committee (project ID 29876).